\documentclass[11pt]{amsart}

\usepackage[a4paper,margin=1.1in]{geometry}
\usepackage{amsmath,amsthm,amssymb,mathtools}
\usepackage{hyperref}
\usepackage{enumitem}
\usepackage [autostyle, english = american]{csquotes}
\MakeOuterQuote{"}

\newtheorem{thm}{Theorem}[section]
\newtheorem{lem}[thm]{Lemma}

\theoremstyle{definition}
\newtheorem{definition}[thm]{Definition}

\numberwithin{equation}{section}

\newcommand{\R}{\mathbb{R}}
\newcommand{\Z}{\mathbb{Z}}

\title{On the best constant in the finitary Vitali covering lemma for high dimensional cubes}
\author{Gian Maria Dall'Ara}
\thanks{Istituto Nazionale di Alta Matematica ``F. Severi" and Scuola Normale Superiore}

\date{\today}

\begin{document}
\maketitle

\begin{abstract} Let $\Gamma_d$ be the largest constant such that every finite collection of cubes in $\R^d$ whose sides are parallel to the coordinate axes admits a disjoint sub-collection occupying a fraction $\Gamma_d$ of its volume. Vitali's greedy algorithm shows that $\Gamma_d\geq 3^{-d}$, and cutting a cube into its $2^d$ dyadic sub-cubes gives $\Gamma_d\leq 2^{-d}$. The question of determining the value of $\Gamma_d$ was first raised by T.~Radó in a 1927 letter to Sierpinski. 
	
	In this paper we investigate the asymptotic behavior of $\Gamma_d$ in the high-dimensional limit. We prove that there exists an absolute constant $c>0$ such that \[
	\Gamma_d\geq c\frac{2^{-d}}{d\log d}
	\] in all dimensions $d$, a significant asymptotic improvement of earlier results by R.~Rado (1949) and Bereg--Dumitrescu--Jiang (2010). This gives an answer to problem D6 in Croft--Falconer--Guy's book "Unsolved problems in geometry". 
\end{abstract}

\section{Introduction}
\label{sec:intro}

In 1908, G.~Vitali \cite{Vitali1908} proved the following well-known "covering lemma", which we state in its finitary version. If $E\subseteq \R^d$ is measurable, we denote by $|E|$ its Lebesgue measure.

\begin{thm}[Vitali covering lemma for axes-parallel cubes \cite{Vitali1908}]\label{thm:vitali}
Any finite collection $\mathcal{C}$ of cubes in $\R^d$ with sides parallel to the coordinate axes admits a disjoint sub-collection $\mathcal{S}$ such that \begin{equation}\label{eq:vitali}
|\cup_{Q\in \mathcal{S}}Q|\geq 3^{-d}|\cup_{Q\in \mathcal{C}}Q|. 
\end{equation}
	\end{thm}

Vitali's proof is based on a greedy algorithm which selects the largest cube (breaking ties arbitrarily), removes from the collection all the cubes intersecting it, and iterates until the collection is empty. This argument applies equally well to a finite collection $\mathcal{C}$ of balls in an arbitrary metric space, showing that $\mathcal{C}$ admits a disjoint sub-collection $\mathcal{S}$ such that \[\cup_{B\in \mathcal{S}}3B\supseteq\cup_{B\in \mathcal{C}}B ,\] where $3B$ is the ball with same center as $B$ and three times the radius. Notice that the latter is a purely metric statement, itself often called the Vitali covering lemma, from which Theorem \ref{thm:vitali} follows easily: one uses the "doubling" property $|3B|=3^d|B|$ and the fact that axes-parallel cubes are the balls of the normed space $\ell_\infty^d$ (that is, $\R^d$ equipped with the norm $\max_{1\leq k\leq d}|x_k|$). In fact, one can prove in the same way an analogue of Theorem \ref{thm:vitali} (with identical constant) in any $d$-dimensional normed space.\newline 

Denote by $\Gamma_d$ the largest constant such that Theorem \ref{thm:vitali} holds with $3^{-d}$ replaced by $\Gamma_d$. Observing that the collection \begin{equation}\label{eq:cell}
\mathcal{C}=\left\{\prod_{j=1}^d [k_j, k_j+1]\colon\, k_j\in \{0,1\}^d\right\}
\end{equation} consists of $2^d$ pairwise intersecting unit cubes, one concludes that  \begin{equation}\label{eq:vitali_interval}
3^{-d}\leq \Gamma_d\leq 2^{-d}.
\end{equation} \newline 

In a 1927 letter to Sierpinski \cite{Radob1928}, T.~Radó observed that $\Gamma_1=\frac{1}{2}$ and posed the problem of determining $\Gamma_2$, expressing his belief that $\Gamma_2=\frac{1}{4}$. In the ensuing decades some progress was made on this problem and its higher-dimensional generalization by Sokolin \cite{Sokolin1940}, R.~Rado \cite{Rado1949I}, Besicovitch (unpublished, reported by R.~Rado), Nordlander \cite{Nordlander1958} and Zalgaller \cite{Zalgaller1960}. Among other things, these authors independently discovered that T.~Radó's conjecture is true in any dimension \emph{if one limits consideration to collections of congruent cubes}. More precisely, let $\gamma_d$ be the maximal constant such that the following holds: every finite collection $\mathcal{C}$ of \emph{congruent} cubes in $\R^d$ with sides parallel to the coordinate axes admits a disjoint sub-collection $\mathcal{S}$ such that \[
|\cup_{Q\in \mathcal{S}}Q|\geq \gamma_d|\cup_{Q\in \mathcal{C}}Q|. 
\] Then we have the following result. 

\begin{thm}[Optimal Vitali covering lemma for axes-parallel congruent cubes]\label{thm:optimal_congruent} $\gamma_d=2^{-d}$.
\end{thm}

It seems plausible that the authors of this theorem viewed it as a special case of the yet-unproven statement that $\Gamma_d=2^{-d}$ (for every $d$), which indeed is called a "long-standing conjecture" in R.~Rado's 1968 paper \cite{Rado1968III}. Remarkably, this conjecture has been disproven in 1973 by Ajtai \cite{Ajtai1973}, who exhibited an explicit finite collection $\mathcal{C}$ of axes-parallel squares in the plane such that every disjoint sub-collection of $\mathcal{C}$ occupies a fraction $<\frac{1}{4}$ of its measure. See the more recent work of Bereg, Dumitrescu and Jiang \cite{Bereg_Dumitrescu_Jiang} for a quantitative improvement of Ajtai's counterexample. \newline 

In contrast to the situation for $\gamma_d$, very little seems to be known about $\Gamma_d$. Slight improvements over Vitali's constant $3^{-d}$ have been established in every dimension by R.~Rado \cite{Rado1949I}, who proved the inequality \[
\Gamma_d \geq (3^d-7^{-d})^{-1}, 
\] and more recently by Bereg, Dumitrescu and Jiang \cite{Bereg_Dumitrescu_Jiang}, who showed that $\Gamma_d\geq 1/\lambda_d$, where $\lambda_d$ is the unique solution in the interval $[(5/2)^d, 3^d]$ of the equation $3^d-(\lambda^\frac{1}{d}-2)^d/2=\lambda$. Notice that $3^d-\lambda_d\in[0,1/2]$, so neither of these results gives a constant $b$ strictly less than $3$ such that $\Gamma_d\geq b^{-d}$ when $d$ is large, thus leaving a huge gap between known lower and upper bounds. In our main result, we show that one may take $b$ as close to $2$ as one pleases. 

\begin{thm}[Almost optimal Vitali covering lemma for high-dimensional axes-parallel cubes]\label{thm:main}
	There exists an absolute constant $c>0$ such that 
	\begin{equation}\label{eq:new}
	\Gamma_d\ \geq c\frac{2^{-d}}{d\log d}\qquad \forall d.
	\end{equation}
	More explicitly, any finite collection $\mathcal{C}$ of cubes in $\R^d$ with sides parallel to the coordinate axes admits a disjoint sub-collection $\mathcal{S}$ such that \begin{equation}\label{eq:main}
		|\cup_{Q\in \mathcal{S}}Q|\geq c\frac{2^{-d}}{d\log d}\, |\cup_{Q\in \mathcal{C}}Q|. 
	\end{equation} The bound continues to hold if one replaces $c$ with $e^{-1-\frac{\log \log d}{\log d}+O\left(\frac{1}{\log d}\right)}$. 
\end{thm}

This gives an answer to Problem D6 in \cite{Croft_Falconer_Guy1991} that, in view of \eqref{eq:vitali_interval}, is almost optimal in the high-dimensional limit. Notice that, as a corollary to the above theorem, we have \[
\lim_{d\rightarrow \infty} \frac{\log(1/\Gamma_d)}{d} = \log 2. 
\] See Section \ref{sec:numerics} for a more refined numerical analysis of the bound given by our argument.

Our proof of Theorem \ref{thm:main} is based on a reduction of the analysis of $\Gamma_d$ to that of $\gamma_d$. This reduction holds in the more general setting of a finite-dimensional normed space, as explained below.

\section{Approximate reduction to unit scale and proof of the main theorem}

In this section, we put our problem in the general context of finite-dimensional normed spaces. A more general framework is discussed in Rado's papers \cite{Rado1949I, Rado1951II, Rado1968III}. \newline 

Let $X$ be a finite-dimensional normed space of dimension $d$. By a ball in $X$ we always mean a closed ball of positive radius. If $B$ is a ball and $\lambda>0$, then $\lambda B$ denotes the ball with same center and $\lambda$ times the radius. We keep denoting the Lebesgue measure of $E\subseteq X$ by $|E|$. The precise normalization of Lebesgue measure is immaterial, as we will only deal with ratios of measures. We have the doubling property $|\lambda B|=\lambda^d|B|$ for every ball $B$. \newline 

Given a finite collection of balls $\mathcal{C}$ in $X$, define 
\[
\Phi(\mathcal{C}):= \max\left\{\frac{|\cup_{B\in \mathcal{S}}B|}{|\cup_{B\in \mathcal{C}}B|}: \mathcal{S}\subseteq \mathcal{C}\ \text{is disjoint}\right\}.
\]
For a set $S\subset(0,\infty)$ of admissible radii, we denote by $\mathcal{A}(X; S)$ the set of all finite collections of balls in $X$ with radii in $S$. Define
\[
\Gamma(X; S):=\inf\bigl\{\Phi(\mathcal{C}) : \mathcal{C}\in\mathcal{A}(X; S)\bigr\}
\] and \[
\Gamma(X):=\Gamma(X; (0,\infty)), \qquad \gamma(X):=\Gamma(X; \{1\}). 
\] 
We remark that, if $B_1$ is the unit ball in $X$, then the latter two constants are denoted $F(B_1)$ and $f(B_1)$ respectively by R.~Rado \cite{Rado1949I}. As mentioned above, we have the following analogue of Theorem \ref{thm:vitali}.

\begin{thm}[Vitali covering lemma on a general normed space]\label{thm:vitali_normed} Let $X$ be a normed space of dimension $d$. Then $\Gamma(X)\geq 3^{-d}$.
	\end{thm} 
Considering all the unit balls containing the origin, one may also easily see that $\gamma(X)\leq 2^{-d}$. Cf.~\cite[Theorem 4]{Rado1949I}. The following properties of the constants $\Gamma(X; S)$ are also easy to verify: \begin{enumerate}
	\item if $S\subseteq T$ then $\Gamma(X; T)\leq \Gamma(X; S)$ (monotonicity), in particular $\Gamma(X)\leq \gamma(X)$; 
	\item $\Gamma(X; tS)=\Gamma(X; S)$ for all $t>0$ (scale-invariance); 
\item $\Gamma_d=\Gamma(\ell_\infty^d)$ and $\gamma_d=\gamma(\ell_\infty^d)$. 
	\end{enumerate} 
The main contribution of this paper is the following theorem. 

\begin{thm}[Approximate reduction to unit scale]\label{thm:many_vs_unit} Let $X$ be a normed space of dimension $d\geq 2$ and let $L$ be a positive integer. Then 
	\begin{equation}\label{eq:many_vs_unit}
		\Gamma(X)\geq (L+2)^{-1}\,\left((2L)^\frac{1}{L+1}(1+L^{-1})\right)^{-d}\, \gamma(X).
	\end{equation} 
In particular, choosing $L=d\log d+O(1)$, one gets 	\begin{equation}\label{eq:many_vs_unit_optimized}
	\Gamma(X)\geq e^{-1-\frac{\log \log d}{\log d}+O\left(\frac{1}{\log d}\right)} \frac{\gamma(X)}{d\log d}. 
	\end{equation} 
\end{thm}

Theorem \ref{thm:main} is obtained plugging the information of Theorem \ref{thm:optimal_congruent} into Theorem \ref{thm:many_vs_unit}, applied to $X=\ell_\infty^d$. \newline 

Fix a normed space $X$ of dimension $d\geq 2$. The proof of Theorem \ref{thm:many_vs_unit} is based on two lemmas. 

\begin{lem}[Comparable versus unit radii]\label{lem:window}
If $r>0$ and $\mu\geq 1$, then
	\[
	\Gamma\big(X; [r,\mu r]\big)\ \geq \mu^{-d}\gamma(X).
	\] 
\end{lem}

\begin{proof}
	Let $\mathcal{C}$ be any finite collection of balls with radii in $[r,\mu r]$. For each $B\in\mathcal{C}$, let $\widetilde B$ be the concentric ball of radius $\mu r$. Let $\widetilde{\mathcal{C}}:=\{\widetilde B\colon\, B\in\mathcal{C}\}$. By definition of $\gamma(X; \{\mu r\})$, there exists a disjoint sub-collection $\widetilde{\mathcal{S}}\subseteq\widetilde{\mathcal{C}}$ with
	\[
	|\cup_{\widetilde B\in \widetilde{\mathcal{S}}} \widetilde B| \geq \gamma(X; \{\mu r\})\, |\cup_{\widetilde B\in \widetilde{\mathcal{C}}}\widetilde B|=\gamma(X)\, |\cup_{\widetilde B\in \widetilde{\mathcal{C}}}\widetilde B|, 
	\]
	where we used the scale-invariance property. 
	
	The collection $\mathcal{S}:=\{B\in \mathcal{C}:\ \widetilde B\in\widetilde{\mathcal{S}}\}$ is clearly disjoint. By the disjointness of $\widetilde{\mathcal{S}}$, 
	\[
	|\cup_{B\in \mathcal{S}}B|
	\geq  \mu^{-d}|\cup_{\widetilde B\in \widetilde{\mathcal{S}}}\widetilde B|
	\geq \mu^{-d}\gamma(X)\, |\cup_{\widetilde B\in \widetilde{\mathcal{C}}}\widetilde B|
	\ \geq \mu^{-d}\gamma(X)\,|\cup_{B\in \mathcal{C}}B|
	\]
	Hence $\Phi(\mathcal{C})\geq  \mu^{-d}\gamma(X)$, and taking the infimum over $\mathcal{C}$ proves the claim.
\end{proof}

For the second lemma, we first need a definition. 

\begin{definition}[(\(\lambda,\mu\))-lacunary sets of radii]
Fix parameters $\lambda\geq 1$ and $\mu\geq 1$. A subset $S\subset(0,\infty)$ is \emph{$(\lambda,\mu)$-lacunary} if there exist intervals $[r_j,s_j]\subseteq (0,\infty)$ ($j=1,2,\ldots$) with
\[
r_{j+1}\ge \lambda\, s_j\quad\text{and}\quad s_j\le \mu r_j \qquad \forall j
\]
such that $S\subset \bigcup_{j\ge1}[r_j,s_j]$.
\end{definition}

\begin{lem}[Lacunary versus unit radii]\label{lem:lacunary}
Let $S\subset(0,\infty)$ be $(\lambda,\mu)$-lacunary, where $\lambda> 1$ and $\mu\geq 1$. Then
\[
\Gamma(X; S)\ \geq \mu^{-d}\, \bigl(1+2\lambda^{-1}\bigr)^{-d}\,  \gamma(X).
\]
\end{lem}

\begin{proof}
Let $\mathcal{C}$ be any finite collection of balls with radii in $S$. By the $(\lambda,\mu)$-lacunarity, $S\subset\bigcup_{j\geq 1} [r_j,s_j]$ with $s_j\leq \mu r_j$ and $r_{j+1}\geq \lambda s_j$. Since $\mathcal{C}$ is finite, there exists a positive integer $m$ such that $S$ intersects $[r_j,s_j]$ only for $j\leq m$. 

Let $\mathcal{C}_j$ be the sub-collection of $\mathcal{C}$ consisting of balls with radii in $[r_j,s_j]$. Define $\mathcal{C}'_m:=\mathcal{C}_m$ and by downward recursion, for $j=m-1,\dots,1$, let $\mathcal{C}'_j$ be the set of balls in $\mathcal{C}_j$ that do not meet any ball in $\bigcup_{k>j}\mathcal{C}'_k$. Set $\mathcal{C}':=\bigcup_{j=1}^m \mathcal{C}'_j$. We have the following two properties:

\smallskip
\noindent\emph{(i)} If $B_1\in\mathcal{C}'_j$ and $B_2\in\mathcal{C}'_k$ with $j\ne k$, then $B_1$ and $B_2$ are disjoint. This is evident by construction. 

\smallskip
\noindent\emph{(ii)} We have 
\[
\cup_{B\in \mathcal{C}} B\ \subseteq\ \cup_{B'\in\mathcal{C}'} (1+2\lambda^{-1})\,B'.
\]
To see this, fix $B\in\mathcal{C}_j\setminus \mathcal{C}'_j$. Then $B$ intersects some $B'\in \mathcal{C}'_k$ with $k>j$. Because $B$ has radius at most $s_j$, $B'$ has radius at least $r_k$, and $r_{k}\geq \lambda s_{k-1}\geq \lambda s_j$, the radius of $B'$ is at least $\lambda$ times that of $B$. Hence $B\subset (1+2\lambda^{-1})B'$.

\smallskip
Applying Lemma~\ref{lem:window} to $\mathcal{C}''_j:=\{(1+2\lambda^{-1})B'\colon\, B'\in\mathcal{C}'_j\}$, which is a collection of balls with radii in $[(1+2\lambda^{-1})r_j,(1+2\lambda^{-1})s_j]\subset [(1+2\lambda^{-1})r_j,\, \mu (1+2\lambda^{-1})r_j]$, gives a disjoint sub-collection $\mathcal{S}''_j\subseteq\mathcal{C}''_j$ with \begin{equation}\label{eq:S''}
|\cup_{B''\in \mathcal{S}''_j} B''|\geq \mu^{-d}\gamma(X)|\cup_{B''\in \mathcal{C}''_j} B''|.
\end{equation}  Summing over $j$, we get \begin{eqnarray}\notag 
\sum_{j=1}^m|\cup_{B''\in \mathcal{S}''_j} B''|&\geq& \mu^{-d}\gamma(X)|\cup_{j=1}^m\cup_{B''\in \mathcal{C}''_j} B''|\\
\notag&=&\mu^{-d}\gamma(X)|\cup_{B'\in \mathcal{C}'} (1+2\lambda^{-1})B'|\\
\label{eq:lac1}&\geq &\mu^{-d}\gamma(X)|\cup_{B\in \mathcal{C}} B|, 
\end{eqnarray}
where in the last inequality we used property (ii) above. 

Let $\mathcal{S}'_j=\{B'\in \mathcal{C}'_j\colon\, (1+2\lambda^{-1})B'\in \mathcal{S}_j''\}$, which is clearly disjoint. By property (i), the collection $\mathcal{S}'=\cup_{j=1}^m\mathcal{S}'_j$ is also disjoint. Thus, \begin{equation}\label{eq:lac2}
|\cup_{B'\in \mathcal{S}'} B'| = \sum_{j=1}^m|\cup_{B'\in \mathcal{S}'_j} B'|= (1+2\lambda^{-1})^{-d}\sum_{j=1}^m|\cup_{B''\in \mathcal{S}''_j} B''|, 
\end{equation}
where we also used the disjointness of each $\mathcal{S}''_j$.
Combining \eqref{eq:lac1} and \eqref{eq:lac2}, we conclude that \[
|\cup_{B'\in \mathcal{S}'} B'| \geq  \mu^{-d}\, (1+2\lambda^{-1})^{-d}\, \gamma(X)|\cup_{B\in \mathcal{C}} B|.
\] Since $\mathcal{S}'$ is a disjoint sub-collection of $\mathcal{C}$, this shows that $\Phi(\mathcal{C})\geq \mu^{-d}\, (1+2\lambda^{-1})^{-d}\, \gamma(X)$ for all $\mathcal{C}\in \mathcal{A}(X; S)$ with $S$ $(\lambda, \mu)$-lacunary, which is the thesis. 
\end{proof}

In order to obtain a proof of Theorem \ref{thm:many_vs_unit}, we now combine Lemmas~\ref{lem:window} and~\ref{lem:lacunary} with a simple pigeonhole argument. The quantity $L$ in the statement of Theorem \ref{thm:many_vs_unit} is essentially the number of pigeonholes used.

\begin{proof}[Proof of Theorem \ref{thm:many_vs_unit}] Fix an integer $J\geq 2$ and a real number $\lambda> 1$. For each residue class $i\in\{0,1,\dots,J-1\}$ define
\[
S_i := \bigcup_{k\in\mathbb{Z}} \bigl[\lambda^{kJ+i},\ \lambda^{kJ+i+1}\bigr].
\]
Clearly, $\bigcup_{i=0}^{J-1}S_i=(0,\infty)$. Let $\mathcal{C}$ be any finite collection of balls. Partition it into the sub-collections \[\mathcal{C}_i:=\{B\in\mathcal{C}:\ B \ \text{has radius in } S_i\}.\] 
Because $|\cup_{B\in \mathcal{C}}B|\leq \sum_{i=0}^{J-1}|\cup_{B\in \mathcal{C}_i}B|$, there exists $i_\ast$ with \[
|\cup_{B\in \mathcal{C}_{i_\ast}}B|\geq J^{-1}|\cup_{B\in \mathcal{C}}B|. 
\]

For a fixed $i$, the set $S_i$ is $(\lambda^{J-1},\lambda)$-lacunary. Applying Lemma~\ref{lem:lacunary} to $\mathcal{C}_{i_\ast}$ yields a disjoint sub-family $\mathcal{S}\subseteq\mathcal{C}_{i_\ast}$ with
\[
|\cup_{B\in \mathcal{S}}B|\ge\ \lambda^{-d}\,\bigl(1+2\lambda^{1-J}\bigr)^{-d}\,\gamma(X)\ |\cup_{B\in \mathcal{C}_{i_\ast}}B|
\ \ge\ J^{-1}\,\lambda^{-d}\,\bigl(1+2\lambda^{1-J}\bigr)^{-d}\,\gamma(X)\ |\cup_{B\in \mathcal{C}}B|.
\]
Since $\mathcal{C}$ is arbitrary, we obtain
\begin{equation}\label{eq:master}
\Gamma(X) \ \ge\ J^{-1}\,\bigl(\lambda(1+2\lambda^{1-J})\bigr)^{-d}\,\gamma(X)
\end{equation} for every $\lambda\geq 1$ and every integer $J\geq 2$. Now we optimize in $\lambda$ for $J$ fixed. Notice that \[
\min_{\lambda\geq1} \lambda(1+2\lambda^{1-J})=\begin{cases}3\qquad &(J=2)\\(2(J-2))^\frac{1}{J-1}\frac{J-1}{J-2}\qquad &(J\geq 3)\end{cases}
\] The estimate obtained for $J=2$ is worse than the Vitali bound of Theorem \ref{thm:vitali_normed}. For $J\geq 3$, after setting $J=L+2$, inequality \eqref{eq:master} yields \eqref{eq:many_vs_unit}. 

Now we consider the regime $L\rightarrow \infty$. A Taylor expansion shows that $(2L)^\frac{1}{L+1}(1+L^{-1})=1+\frac{\log L}{L}+O(L^{-1})$. Let $L=d\log d+O(1)$ ($L$ is an integer, so the error term is needed). A simple calculation gives \[
(2L)^\frac{1}{L+1}(1+L^{-1})=1+\frac{1}{d}+\frac{\log \log d}{d\log d}+O\left(\frac{1}{d\log d}\right)
\]
and  \[
(L+2)\left((2L)^\frac{1}{L+1}(1+L^{-1})\right)^d=d\log d \exp\left(1+\frac{\log \log d}{\log d}+O\left(\frac{1}{\log d}\right)\right).
\]
Plugging this into \eqref{eq:many_vs_unit} we obtain \eqref{eq:many_vs_unit_optimized}. 
\end{proof}

\section{Numerics}\label{sec:numerics}

In a fixed dimension $d$, one may numerically maximize the RHS of \eqref{eq:many_vs_unit}, obtaining a result that is more precise than estimate \eqref{eq:new}. In this section we perform an elementary numerical analysis of the constants in Theorem \ref{thm:many_vs_unit} that shows in particular that \emph{our argument improves the classical Vitali covering lemma for cubes in every dimension $d\geq 14$}. 

Define \[
g(x)=(2x)^\frac{1}{x+1}(1+x^{-1}),\qquad h_d(x)=(x+2)g(x)^d\qquad (x>0). 
\] Since $\lim_{x\rightarrow \infty}h_d(x)=\infty$, there exists a minimal positive integer $L_d$ such that $h_d(L_d)\leq h_d(L)$ for all $L\in \Z^+$. By Theorem \ref{thm:optimal_congruent} and Theorem \ref{thm:many_vs_unit}, we have the inequality \begin{equation}\label{eq:best_bound}
\Gamma_d\geq 2^{-d}h_d(L_d)^{-1}.
\end{equation}

The task of computing numerically the quantities $L_d$ and $h_d(L_d)$ for a fixed $d$ is made easier by the crude bound \begin{equation}\label{eq:L_d_bound}
L_d\leq 4d\log (4d)+1. 
\end{equation}
\begin{proof}[Proof of \eqref{eq:L_d_bound}]
We have to check that the logarithmic derivative \[\frac{h_d'(x)}{h_d(x)}=\frac{1}{x+2}-d\frac{\log(2x)}{(x+1)^2}\] is positive for $x>4d\log (4d)$. 
By the formula, we see that the logarithmic derivative is positive on the set $P=\{x>0\colon\, (x+1)^2>d(x+2)\log(2x)\}$. Using the bounds $\log(2x)\leq 2\log x$ ($x\geq 2$) and $(x+1)^2\geq (x+2)^2-2(x+2)$, we see that $P\supseteq P_1=\{x>\max\{2d\log x, 2\}\}$. Finally, using the change of variable $x=t\log t$ and the estimate $\log\log t\leq \log t$ ($t>1$), we obtain $P_1\supseteq P_2=\{t>4d\}=\{x>4d\log (4d)\}$. \end{proof}

By \eqref{eq:L_d_bound}, we may compute $L_d$ and $m_d=2^dh_d(L_d)$ for $d\leq 20$ by evaluating $h_d(L)$ for $L\leq 4d\log(4d)+1$. The following table was obtained using Mathematica. 

\begin{center}
\begin{tabular}{|r| r| r| r|}
	\hline
	$d$ & $L_d$ & $m_d$ & $m_d/3^d$\\
	\hline
	1 & 1 & 16.971 & 5.657\\
	2 & 4 & 86.152 & 9.572\\
	3 & 8 & 287.026 & 10.631\\
	4 & 13 & 818.895 & 10.110\\
	5 & 18 & 2153.470 & 8.862\\
	6 & 23 & 5379.103 & 7.379\\
	7 & 28 & 12971.417 & 5.931\\
	8 & 34 & 30486.998 & 4.647\\
	9 & 39 & 70264.112 & 3.570\\
	10 & 45 & 159472.691 & 2.701\\
	11 & 51 & 357492.434 & 2.018\\
	12 & 57 & 793261.059 & 1.493\\
	13 & 63 & 1745233.682 & 1.095\\
	\hline 
	14 & 69 & 3811881.534 & 0.797\\
	15 & 75 & 8274033.774 & 0.577\\
	16 & 81 & 17862582.185 & 0.415\\
	17 & 88 & 38379142.946 & 0.297\\
	18 & 94 & 82115993.482 & 0.212\\
	19 & 101 & 175038063.601 & 0.151\\
	20 & 107 & 371863945.976 & 0.107\\
	\hline
\end{tabular}
\end{center}
As highlighted in the table, our estimate is worse than Vitali's when $d\leq 13$, while it improves it when $14\leq d\leq 20$. 
We now conclude our numerical analysis by showing that estimate \eqref{eq:best_bound} is better than Vitali's bound in all dimension $d\geq 14$, that is, \begin{equation}\label{eq:improvement} m_d< 1\end{equation} for all $d\geq 14$. Assume that \eqref{eq:improvement} holds in dimension $d_0$ and that \begin{equation}\label{eq:g_bound} g(L_{d_0})\leq \frac{3}{2}.\end{equation} Then for $d\geq d_0$ we have \begin{eqnarray*}
	m_d&=&\left(\frac{2}{3}\right)^dh_d(L_d)	\\
&\leq& \left(\frac{2}{3}\right)^dh_d(L_{d_0})\\
	&=&\left(\frac{2}{3}\right)^d(L_{d_0}+2)g(L_{d_0})^d\\
	&=& m_{d_0}\left(\frac{2}{3}\right)^{d-d_0}g(L_{d_0})^{d-d_0}<1.
\end{eqnarray*} Thus, if we get an improvement over Vitali's bound in dimension $d_0$ and \eqref{eq:g_bound} holds, then we also get an improvement in all higher dimensions. 

The function $g$ is strictly decreasing on $[1,\infty)$, as may be easily seen by a logarithmic derivative computation (while this is not needed for our purposes, we take the opportunity to point out that the monotonicity of $g$ also implies that $L_d$ is non-decreasing in $d$
). Since one may numerically verify that $g(9)<\frac{3}{2}<g(8)$, condition \eqref{eq:g_bound} is equivalent to $L_d\geq 9$. 
Since $L_{14}=69$ (as computed above), our proof that $m_d<1$ in all dimensions $d\geq 14$ is complete.


\end{document}